\documentclass[12pt,a4paper,draft]{article}
\usepackage{amsmath,amssymb,amsfonts, amsthm}
\date{\begin{flushleft} {\small } \end{flushleft}}
\usepackage[cp1251]{inputenc}
\usepackage[english,russian]{babel}
\usepackage{amssymb,amsmath}

\newcounter{sbpr}

\newtheorem{thm}{Theorem}

\newtheorem{utt}{Proposition}
\newtheorem{cor}{Corollary}
\newtheorem{lm}{Lemma}
\theoremstyle{remark}
\newtheorem{rem}{Remark}
\theoremstyle{definition}
\newtheorem{df}{Definition}

\newtheorem{ex}{Example}
\newtheorem{qu}{Question}

 \def\proof{{\textbf{Proof}. }}
 \def\T{{\mathbb T}}
\def\C{{\mathbb C}}

\def\R{{\mathbb R}}

\def\codim{{\rm codim\:}}
\def\codimc{{\rm codim}_\C\:}

\def\1{{1}_\T}

\def\({{\rm (}}
\def\){{\rm )}}

\def\mdskip{\vskip-\lastskip\vskip\medskipamount}

\def\QED{{\parfillskip0pt\hfil$\square$\par}\mdskip}
\title{On piecewise pluriharmonic functions\thanks {Supported by NSh-4850.2012.1} }
\author{B. Kazarnovskii}
\begin{document}
\maketitle {\footnotesize
We extend some results on piecewise linear functions on $\C^n$
to piecewise pluriharmonic functions on any complex manifold.
 We construct a ring
 generated by currents $h$ and $dd^ch$,
 where $\{h\}$ is a finite set of piecewise pluriharmonic functions.
 We prove
 that, with some restrictions on the set $\{h\}$, the map
 $\{h\mapsto dd^ch,\ dd^ch\mapsto0\}$
 can be continued to
 the derivation on the ring.
 As a corollary,
 the current $dd^cg_1\wedge\cdots\wedge dd^cg_k$
 depends on the product of piecewise pluriharmonic functions $g_1,\cdots,g_k$  only.
}
 \section{Results. \label{1}}
  A function $g\colon M\to\R$ on a complex manifold  $M$ is called pluriharmonic if $dd^cg=0$
  (recall that $d^cg(x_t)=dg(\sqrt{-1}x_t\)$).
  Another definition: the function g is a real part of some holomorphic function in some neighborhood of any point $x\in M$.
   \begin{df}
   Continuous function $g\colon M\to\R$ on $n$-dimensional complex manifold  $M$ is called \emph{piecewise pluriharmonic}
   (or PPH-function)
   if $g$ is pluriharmonic on any closed $2n$-dimensional
   simplex of some locally finite triangulation of the manifold $M$
   (see subsection 2.1).
    \label{dfpph}
 \end{df}
Piecewise linear functions on the space $\C^n$ are the simplest examples of PPH-functions.
Piecewise linear functions  are used in
convex geometry, algebraic geometry, and complex analysis \cite{monge}-\cite{Pol}.
We extend some results on piecewise linear functions on $\C^n$ \cite{monge}
to PPH-functions on any complex manifold.

There exists a nonconstant PPH-function on any complex manifold.
Indeed,
let $h$ be a piecewise linear function on the space $\C^n$
and $h=0$ outside some small neighborhood of zero.
Let $H$ be a function on $M$ such that $H=h$ on the coordinate neighborhood of some point $x\in M$
and $H=0$ outside this neighborhood.
The function $H$ is piecewise pluriharmonic.
   \begin{df}
   Let $A$ be some vector space
   consisting of PPH-functions
   on an $n$-dimensional complex manifold $M$.
   The vector space $A$ is called
   \emph{the constructive space} if for any finite set of elements $H=\{h_i\in A\}$
   there exists a triangulation ${\cal P}_H$  of $M$ such
   that each function $h_i$ is pluriharmonic on every closed $2n$-dimensional
  simplex $\Delta\in{\cal P}_H$.
  \label{dfconstr}
 \end{df}
 In the sequel, all PPH-functions are the elements of some fixed constructive space $A$.
   \begin{ex}
  The space of piecewise linear functions on $\C^n$ is the constructive space.
  \label{exconstr}
 \end{ex}
  \begin{qu}
  Is the space of all PPH-functions constructive?
  \label{quconstr}
 \end{qu}
Let $h_i\colon M\to\R$  be a function on a complex $n$-dimensional manifold $M$.
The mixed Monge-Ampere operator on $M$ of degree $k$ is (by definition)
the map $(h_1,\cdots,h_k)\mapsto dd^ch_1\wedge\cdots dd^ch_k$.
If $h_1,\cdots,h_k$ are continuous plurisubharmonic functions,
 then \cite{BT} the Monge-Ampere operator value $dd^ch_1\wedge\cdots\wedge dd^ch_k$
  is well defined as a current
  (that is a functional on the space of smooth compactly supported differential $(2n-2k)$-forms).
  This means that if the sequence of smooth plurisubharmonic functions $(f_i)_j$ converges locally uniformly to $h_i$,
  then the sequence of currents $dd^c(f_1)_j\wedge\cdots dd^c(f_k)_j$
  converges to the limit current,
  independent of the choice of approximation.
  This limit current is the current of measure type,
  ie may be continued to a functional on the space of continuous compactly supported forms.
  It follows that any polynomial in the variables $h_1,\cdots,h_p,dd^cg_1,\cdots,dd^cg_q$
  with continuous plurisubharmonic functions $h_i$ and $g_i$ gives
  the well defined current on $M$.

  It is easy to prove that any PPH-function can be locally written as a difference of two
  continuous plurisubharmonic functions.
  It follows that above defined currents are well defined for PPH-functions $h_1,\cdots,h_p,g_1,\cdots,g_q$ also.
  %
  \begin{thm}
  Let $\cal A$ be a ring of currents,
   generated by currents $h$ and $dd^c h$ with $h\in A$.
  There exists a derivation $\delta$ on the ring
  ${\cal A}$ such that $\delta(h)=dd^ch$ and $\delta(dd^ch)=0$
  for any {\rm PPH}-function $h\in A$.
  \label{thmMonge}
 \end{thm}
 \begin{rem}
 Let $h_i$ be PPH-functions and $F=h_0dd^ch_1\wedge\cdots\wedge dd^ch_k$;
 then $\delta(F)=dd^cF$.
 But $\delta(h^2)=0$ and $dd^c(h^2)\not=0$ for any nonconstant pluriharmonic function $h$.
 \end{rem}
 \begin{cor}
  $\delta^k(h_1\cdots h_k)=k!dd^ch_1\wedge\cdots\wedge dd^ch_k$ for any {\rm PPH}-functions
  $h_1,\cdots,h_k$.
 \label{corEst}
 \end{cor}
 \begin{cor}
 The current $dd^ch_1\wedge\cdots\wedge dd^ch_k$
 depends on the product of {\rm PPH}-functions $h_1,\cdots,h_k$  only.
 \label{corMonge}
 \end{cor}
  \begin{cor}
  Let the space $A$,
  generated by the {\rm PPH}-functions
  $h_1,\cdots,h_k$,
  be a constructive space
  and let $h=\max(h_1,\cdots,h_k)$.
  Suppose that $h\in A$. Then $dd^c(h-h_1)\wedge\cdots dd^c(h-h_k)=0$.
 \label{corGusev}
 \end{cor}
 The proof of theorem \ref{thmMonge} uses (as the proofs of similar theorems for piecewise linear functions in \cite{Est,monge})
 a construction of \emph{$k$-th corner locus}.
 Such corner loci (in some more special situation) were constructed in \cite{Est}.
  \section{ PPH-cycles and corner loci. \label{2}}
   \subsection{P-cycles and operator $D_c$. \label{21}}
   \emph{A simplex $\Delta$ in a complex $n$-dimensional manifold} $M$ is the image
  of a smooth nonsingular embedding $\Delta^{2n}\to M$,
  where  $\Delta^{2n}$ is the standard closed $2n$-dimensional simplex.
  A locally finite set $\cal P$ of simplices in $M$
  is called \emph{a triangulation}
  if $\cup_{\Delta\in\cal P}\Delta=M$
  and intersection of any two simplices from $\cal P$
  is their common (may be empty) face.
  Any face of any simplex $\Delta\in\cal P$ is called \emph{a cell} of triangulation $\cal P$.
  An odd form $\omega$ on the cell $\Delta$ is called the \emph{frame of $\Delta$}
  (recall that the odd form on a manifold with orientations $\alpha,\beta$
  is a pair of forms $\omega_\alpha,\omega_\beta$ such that $\omega_\alpha=-\omega_\beta$).
  Let \emph{$k$-dimensional chain} (or $k$-chain)
  be a map $X\colon\Delta\mapsto X_\Delta$ of the set of all $k$-dimensional cells to their frames.
    Let the map $\partial$ take each $k$-chain $X$ to $(k-1)$-chain $\partial X$,
  where
       $$(\partial X)_\Lambda=\sum_{\Delta\supset\Lambda,\dim\Delta=k} X_\Delta,$$
 where the orientations of the cells $\Lambda$ and $\Delta$ agreed as usual.
 A $k$-chain $X$ is called a \emph{$k$-cycle} if $\partial X=0$.
 \begin{cor}
 For any $k$-chain $X$ the $(k-1)$-chain $\partial X$ is a cycle.
 \label{cordsqu}
 \end{cor}
   \emph{In the sequel we assume that $k\geq n$.}
 
  Let $\Delta,T_\Delta^x$, and $\C_\Delta^x$ be (respectively) a $k$-dimensional cell,
  the tangent space of $\Delta$ at the point $x\in\Delta$, and
  the maximal complex subspace of $T_\Delta^x$.
  Say that the point $x\in\Delta$ is nondegenerate
  if $\codimc\C_\Delta^x =\codim T_\Delta^x$.
  For $k=2n,2n-1$ any point of $k$-dimensional cell is nondegenerate.
  If $x$ is nondegenerate, then $\dim T_\Delta^x - \dim_\R \C_\Delta^x = 2n-k$.
  If $x$ is degenerate, then
  the form $X_\Delta$ (from definition \ref{dftpol}) is zero at $x$.
 \begin{df}
  A $k$-cycle $X$ is called \emph{$P$-cycle} if the following conditions hold:

 {\rm(1)}
 $\deg X_\Delta=2n-k$

 {\rm(2)}
 $X_\Delta=\sum_{1\leq j\ll\infty} g^j_{\Delta} W^j_{\Delta}$,
 where  $(2n-k)$-forms $W^j_{\Delta}$ are closed

 {\rm(3)}
  $W_\Delta(\xi_1,\cdots,\xi_{2n-k})=0$
 for any $x\in\Delta$ and any $(\xi_1,\cdots,\xi_{2n-k})\subset T_\Delta^x$
 if $\exists i\colon \xi_i\in\C_\Delta^x$.
 \label{dftpol}
 \end{df}
 Let $X$ be a $k$-dimensional P-cycle.
  \emph{Define the $(k-1)$-cycle $D_cX=\partial Y$,}
 where $Y$ is a $k$-chain such that
 \begin{equation}
 \label{Y}
    Y_\Delta=\sum_{1\leq j\ll\infty} d^c G^j_{\Delta}\wedge W^j_\Delta,
  \end{equation}
  where $G^j_\Delta$ are smooth functions in some neighborhood of $\Delta$
 such that $G^j_\Delta(x)=g^j_{\Delta}(x)$ for any $x\in\Delta$.
 %
 \begin{cor}
 If $k=n$ then $D_cX=0$.
 \label{cormorethann}
 \end{cor}
 \begin{utt}
   {\rm(a)} The $(k-1)$-cycle $D_cX$ does not depend on the choice of functions $G^j_\Delta$ in {\rm(\ref{Y})}
   and on the choice of functions $g^j_\Delta$ in
   decomposition $X_\Delta=\sum_{1\leq j\ll\infty} g^j_{\Delta} W^j_{\Delta}$ from definition {\rm\ref{dftpol}}.

 {\rm(b)}
  $D_cX_{\Lambda}(\xi_1,\cdots,\xi_{2n-k+1})=0$
 for $x\in\Lambda$ and $(\xi_1,\cdots,\xi_{2n-k+1})\subset T_\Lambda^x$
 if $\exists i\colon \xi_i\in\C_\Lambda^x$.
  \label{uttdelta}
 \end{utt}
 \proof
 If $F(x)=0$ for any $x\in\Delta$, then the form $d^cF$ is zero on the subspace $\C^x_\Delta\subset T^x_\Delta$ .
 Using definition \ref{dftpol} (3) we get that the form $d^c F\wedge W^j_\Delta$ is zero on $\Delta$.
 It follows that the form $(D_c X)_\Lambda$
 does not depend on the choice of functions $G^j_\Delta$.

 Now we prove that the form $(D_c X)_\Lambda$
 does not depend on the choice of functions $g^j_\Delta$.
 Let
 $x_1\wedge\cdots\wedge x_{2n-k+1}\not=0$,
 where  $x_i\in T^x_\Delta$ and $x_1\in\C^x_\Delta$.
 If $(\sum_j g^j_{\Delta} W^j_{\Delta})=0$,
 then
  \begin{multline*}
  \sum_j d^c G^j_{\Delta}\wedge W^j_\Delta(x_1,x_2,\cdots)=
   \sum_j dG^j_{\Delta}\wedge W^j_\Delta(ix_1,x_2,\cdots)=\\
   \sum_j dg^j_{\Delta}\wedge W^j_\Delta(ix_1,x_2,\cdots)=
   d\left(\sum_j g^j_{\Delta}\wedge W^j_\Delta\right)(ix_1,x_2,\cdots)=0.
 \end{multline*}
 The first equality follows from the definition of operator $d^c$
 and from definition \ref{dftpol},
 the third -- from the closedness of forms $W^j_\Delta$.
 Assertion (a) is proved.

 Let $\dim\Lambda=k-1$ and let
 $x_1\wedge\cdots\wedge x_{2n-k+1}\not=0$,
 where  $x_i\in T^x_\Lambda$ and $x_1\in\C^x_\Lambda$.
 We prove that $(D_c X)_\Lambda(x_1,\cdots,x_{2n-k+1})=0$.
  \begin{multline*}
 (D_c X)_\Lambda(x_1,x_2,\cdots)=
 \sum_{\Delta\supset\Lambda,\dim\Delta=k} \sum_j d^c G^j_{\Delta}\wedge W^j_\Delta(x_1,x_2,\cdots)=
  \\
 \sum_{\Delta\supset\Lambda,\dim\Delta=k} \sum_j dG^j_{\Delta}\wedge W^j_\Delta(i x_1,x_2,\cdots)=
 \sum_{\Delta\supset\Lambda,\dim\Delta=k}d X_\Delta(i x_1,x_2,\cdots)=
 \\
 d\left(\sum_{\Delta\supset\Lambda,\dim\Delta=k}X_\Delta(i x_1,x_2,\cdots)\right)=
 d(\partial X)_\Lambda=0.
 \end{multline*}
 The first equality is a definition of $(k-1)$-chain $\partial Y$.
 Other equalities
 follow from the definition of operator $d^c$,
the closedness of forms $W^j_\Delta$, and the closedness of the $k$-chain $X$.
\QED
  \begin{rem}
 $D_cX$ is not necessarily a P-cycle.
 \end{rem}
 \subsection{Corner loci of PPH-polynomials. \label{22}}
 \begin{df}
 Let $H=\{h_1,\cdots,h_q\}$ be a basis of the constructive space $A$,
 $P(x_1,\cdots,x_q)$ be a polynomial of degree $m$ in the variables $x_1,\cdots,x_q$.
 The function $P(h_1,\cdots,h_q)$
 is called a {\rm PPH}-polynomial of degree $m$.
 \label{dfpphpol}
 \end{df}
  The degree of PPH-polynomial is not uniquely defined.
  But PPH-polynomials of degree $0$ are constants
  and PPH-polynomials of degree $1$ are PPH-functions.

  Below we fix the set $H$
  and the triangulation ${\cal P}_H$.
  The restriction of PPH-polynomial $P$ to standardly oriented $2n$-cells of triangulation ${\cal P}_H$
   give the $2n$-dimensional P-cycle $X^P$.
  Say that $(2n-1)$-cycle $D_cX^P$ is
   \emph{the corner locus of {\rm PPH}-polynomial} $P$.
    \begin{ex} \emph{Corner locus of {\rm PPH}-function.}
    Any $(2n-1)$-dimensional cell $\Delta$ of triangulation ${\cal P}_H$
    is a common face
    of $2n$--dimensional cells $B^+$ and $B^-$.
    Let $h\in A$ and $h_\Delta^+=h|_{B^+},\,h_\Delta^-=h|_{B^-}$.
    The ordering of the pair $(B^+,B^-)$ sets the coorientation of $\Delta$.
    The standard orientation of $M$ and the coorientation of $\Delta$ together set the orientation of the cell $\Delta$.
    Using this orientation put $(D_cX^h)_\Delta=d^ch_\Delta^+-d^ch_\Delta^-$.
    \label{exHyper}
     \end{ex}
  \begin{cor}
  For any {\rm PPH}-polynomial $P$
 \begin{equation}
 \label{delta1}
    (D_cX^P)_\Delta=\sum_{1\leq i\leq q} \frac {\partial P}{\partial x_i}(h_1,\cdots,h_q) (d^c(h_i)_\Delta^+-d^c(h_i)_\Delta^-),
 \end{equation}
 where $1$-forms $d^c(h_i)_\Delta^\pm$ on the $(2n-1)$-dimensional oriented cell $\Delta$ are defined in the text of example {\rm\ref{exHyper}}.
  \label{cordelta1}
  \end{cor}
 \begin{cor}
 The corner locus $D_cX^P$ is a {\rm P}-cycle.
 \label{cordelta2}
 \end{cor}
 \begin{df}
 The {\rm P}-cycle $D_c^kX^P$  is called a $k$-th corner locus of {\rm PPH}-polynomial $P$.
 \label{dfcornerloc}
 \end{df}
 \emph{The validity of definition \ref{dfcornerloc} is based on the assertion (1) of lemma \ref{lmcorrcor}}.
 Below we use the following notation:
\begin{enumerate}
\item
$I=\{i_1,\cdots,i_q \colon i_j\geq0\}$, $|I|=i_1+\cdots+i_q$, $I!=i_1!\cdots i_q!$,
$x^I=x_1^{i_1}\cdots x_q^{i_q}$.
\item
If $i_p >0$ then $I\setminus p=\{i_1,\cdots,i_p-1,\cdots,i_q\}$;
else $I\setminus p=\emptyset$.
 \item
 If $I\not=\emptyset$ then
$P_I(h_1,\cdots,h_q)=\frac{\partial^{|I|} P}{\partial x_1^{i_1} \cdots \partial x_q^{i_q}}(h_1,\cdots,h_q)$;
else $P_I(h_1,\cdots,h_q)=0$.
\item
$\Delta\mapsto\tilde\Delta$ is some fixed mapping of the set of cells of triangulation ${\cal P}_H$ into itself such
that

(a)  $\Delta\subset\tilde\Delta$

(b) $\dim\tilde\Delta=2n$

(c) if $\dim\Delta=2n$ then $\tilde\Delta=\Delta$.
\item
$H^i_\Delta$ is the restriction of pluriharmonic function $(h_i)_{\tilde\Delta}$
to some neighborhood of the cell $\Delta$
\item
$G^i_{\Delta,\Gamma}$ is the restriction of pluriharmonic function $(h_i)_{\Gamma}$
to some neighborhood of the cell $\Delta$,
where $\Gamma$ ranges over the set of $2n$-dimensional cells
containing $\Delta$.
\end{enumerate}
 \begin{lm}
  For any {\rm PPH}-polynomial $P$ and any $k\geq1$

  {\rm(1)} $D_c^kX^P$ is a P-cycle.

  {\rm(2)} If $\dim\Delta=2n-k$ then
  $$(D_c^{k}X^P)_\Delta=\sum_{|I|=k} P_I(h_1,\cdots,h_q) Q_{I,\Delta}^k,$$
 where $Q_{I,\Delta}^k$ is a polynomial of degree $k$ in variables $\{d^cG^i_{\Delta,\Gamma}\}$
 and is independent of the choice of a polynomial $P$.

 {\rm(3)}
 If $|I|=k$ and $P=x^I$ then $(D_c^{k}X^P)_\Delta=I! Q_{I,\Delta}^k$.
  \label{lmcorrcor}
   \end{lm}
\proof
The proof is by induction on k.
For $k=1$, all the assertions follow from corollary \ref{cordelta1}.
If the assertion is true for $k-1$ then
  \begin{multline*}
 (D_c^{k}X^P)_\Lambda=
 \sum_{\Delta\supset\Lambda,\dim\Delta=2n-k+1}\sum_{|I|=k-1}d^cP_I(H^1_\Delta,\cdots,H^q_\Delta)Q_{I,\Delta}^{k-1}=\\
 \sum_{\Delta\supset\Lambda,\dim\Delta=2n-k+1}\sum_{|I|=k-1}\left(\sum_{1\leq i\leq q}
\frac{\partial P_I}{\partial x_i}(h_1,\cdots,h_q) d^cH^i_\Delta\right)
 \wedge Q_{I,\Delta}^{k-1}=\\
  \sum_{\Delta\supset\Lambda,\dim\Delta=2n-k+1}\sum_{|I|=k}P_I(h_1,\cdots,h_q)
  \sum_{1\leq i\leq q} d^cH^i_\Delta\wedge Q_{I\setminus i,\Delta}^{k-1}=\\
\sum_{|I|=k} P_I(h_1,\cdots,h_q)\sum_{\Delta\supset\Lambda,1\leq i\leq q}d^cH^i_\Delta\wedge Q_{I\setminus i,\Delta}^{k-1},
   \end{multline*}
where the product $d^cH^i_\Delta\wedge Q_{I\setminus i,\Delta}^{k-1}$ is an odd form (as a product of even and odd forms)
restricted to the cell $\Lambda$.
So we can put
$$Q_{I,\Lambda}^k=\sum_{\Delta\supset\Lambda,1\leq i\leq q}d^cH^i_\Delta\wedge Q_{I\setminus i,\Delta}^{k-1}.$$
Assertion (2) is proved.

Applying assertion (2) to $P=x^I$ we get assertion (3).

Combining assertion (3) and proposition \ref{uttdelta} (b) we obtain assertion (1).
 \begin{cor}
 If $k\geq n$ then $D_c^{k}X^P=0$.
 \label{cormorethann2}
 \end{cor}
 \subsection{Corner loci of PPH-cycles. \label{23}}
  Let $i=1,\cdots,p$ and $X_i$ be a $k$-dimensional P-cycle
  (definition \ref{dftpol}).
  For PPH-polynomials $P^1,\cdots,P^p$ we define the $k$-cycle $X=P^1X_1+\cdots+P^pX_p$
  as $X_\Delta=P^1(X_1)_\Delta +\cdots+P^p(X_p)_\Delta$.
 \begin{df}
 The $k$-cycle $X$ is called a {\rm PPH}-cycle if the forms $(X_i)_\Delta$ are closed.
 Put $\deg X=\max_i \deg P^i$.
 \label{dfphp}
 \end{df}
 \begin{cor}
 Any {\rm PPH}-cycle is a P-cycle.
 \label{corPPHT}
 \end{cor}
 \begin{cor}
 {\rm PPH}-cycles form the module over the ring of {\rm PPH}-polynomials.
 \label{cormod}
 \end{cor}
 \begin{df}
 The cycle $D_cX$ is called a corner locus of {\rm PPH}-cycle $X$.
 \label{dfcuclecorner}
 \end{df}
 \begin{utt}
  The corner locus $D_cX$ of any $k$-dumensional {\rm PPH}-cycle $X$ is a $(k-1)$-dimensional {\rm PPH}-cycle  and
  $\deg D_cX=\deg X-1$.
 \label{uttPPHDC}
 \end{utt}
 \proof
Let $Y$ be a PPH-cycle of degree $0$.
Using the definition of operator $D_c$, we have
$$(D_c(h_jY))_\Lambda=\sum_{\Delta\supset\Lambda,\dim\Delta=k} d^cH^j_\Delta\wedge Y_\Delta,$$
where the functions $H^i_\Delta$  defined in subsection \ref{22}.
The forms $d^cH^j_\Delta\wedge Y_\Delta$ are closed.
It follows that $D_c(h_jY)$ is PHP-cycle of degree $0$.
Now it remains to observe that
if $X=P^1X_1+\cdots+P^pX_p$ then
$$D_cX = \sum_{1\leq i\leq p,\:1\leq j\leq q}\frac{\partial P^i}{\partial x_j}(h_1,\cdots,h_q)\:D_c(h_jX_i).$$
\begin{cor}
 If $h\in A$ then $D_c(h^kX)=kh^{k-1}D_c(hX)$.
 \label{corhk}
 \end{cor}
 \addtocounter{sbpr}{1}
 \subsection{PPH-cycles as currents. \label{24}}
Let $X$ be a $k$-chain.
Now suppose $\bar X$ is a current such that
$$\bar X(\varphi)=\sum_{\Delta\in{\cal P}_H,\dim\Delta=k}\int_\Delta X_\Delta\wedge\varphi.$$
 \begin{lm}
  Let $X$ be a $k$-cycle such that the forms $X_\Delta$ are closed.
  Then the current $\bar X$ is closed.
  \label{lmforStokes}
 \end{lm}
\proof
  \begin{multline*}
d\bar X(\psi)=
\bar X(d\psi)=
\sum_{\Delta\in{\cal P}_H,\dim\Delta=k} \int_\Delta X_\Delta\wedge d\psi=
\sum_{\Delta\in{\cal P}_H,\dim\Delta=k} \int_\Delta d(X_\Delta\wedge\psi)=
\\
\sum_{\Lambda\in{\cal P}_H,\dim\Lambda=k-1} \sum_{\Delta\supset\Lambda,\dim\Delta=k} \int_\Lambda X_\Delta\wedge\psi=
\sum_{\Lambda} \int_\Lambda (\partial X)_\Lambda\wedge\psi=
(\overline{\partial X})(\psi) =0
\end{multline*}
%
 \begin{cor}
  Let $X$ be a $k$-dimensional {\rm PPH}-cycle of degree $0$.
  Then the current $\bar X$ is closed.
  \label{corforStokes}
 \end{cor}
  \begin{utt}
  Let $h\in A$, $X$ be a $k$-dimensional {\rm PPH}-cycle of degree $0$.
  Then  $\overline{D_c(hX)}=dd^c(h\bar X)$.
  \label{uttEst2}
 \end{utt}
\proof
  \begin{multline}
   dd^c(h\bar X)(\psi)=
 -\sum_{\Delta\in{\cal P}_H,\dim\Delta=k}\int_\Delta h X_\Delta\wedge dd^c\psi=
 \\
 (-1)^{k+1}\sum_{\Delta}\int_\Delta \left(d(h X_\Delta\wedge d^c\psi) - dh\wedge X_\Delta\wedge d^c\psi\right)
 =\\
 (-1)^{k}\sum_{\Delta}\int_\Delta dh\wedge X_\Delta\wedge d^c\psi
=
(-1)^{k+1}\sum_{\Delta}\int_\Delta d^ch\wedge X_\Delta\wedge d\psi=\\
 \sum_\Delta\int_\Delta d(d^ch\wedge X_\Delta\wedge \psi)=
 \sum_{\Delta}\sum_{\Lambda\subset\Delta,\dim\Lambda=k-1}\int_\Lambda d^ch\wedge
  X_\Delta\wedge \psi=\overline{D_c(hX)}(\psi).
  \label{formuladdc}
 \end{multline}
 The formula (\ref{formuladdc}) consists of seven equalities.
  We shall comment to each one.
 \begin{enumerate}
 \item
  Determination of the current derivative
  and the identity $dd^c=-d^cd$.
\item
  The closedness of the form $X_\Delta$.
\item
  The Stokes formula and the $k$-chain $hX$ closedness.
\item
  Let the form $\psi$ bidegree is $(k-n-1,k-n-1)$
 (the values of the current $dd^c(h\bar X)$ on homogeneous components of other degrees are zero).
 Suppose $x\in\Delta$ and $\xi_1,\cdots,\xi_{2k-2n}$ is a basis of the (real) vector space $\C^x_\Delta$.
 It is easy to prove that $$(dh\wedge d^c\psi)(\xi_1,\cdots,\xi_{2k-2n})=-(d^ch\wedge d\psi)(\xi_1,\cdots,\xi_{2k-2n}).$$
 Hence $(dh\wedge X_\Delta\wedge d^c\psi)=-(d^ch\wedge X_\Delta\wedge d\psi)$.

\item
 The closedness of the forms $d^ch$ and $X_\Delta$.
\item
  The Stokes formula.
\item
  The deternination of $D_c$.
 \end{enumerate}
 \section{Theorem \ref{thmMonge} (proof). \label{3}}
  Below we use the following notation (with the notation of subsection \ref{22}):
\begin{enumerate}
\item
  ${\cal B}$ is the ring of PPH-polynomials.
\item
  ${\cal B}'$ is the symmetric algebra of the space $A$.
\item
${\cal T}$ is the ring (with the unity) generated by currents $dd^ch$,
where $h\in A$.
\item
  ${\cal T}'={\cal T}''/I$, where ${\cal T}''$ is the symmetric algebra of the space of currents $dd^ch$,
 where $h\in A$, and $I$ is the ideal generated by elements of degree $>n$.
\item
  ${\cal A}'={\cal B}'\otimes{\cal T}'$ is the tensor product of the rings.
\item
  $\delta'$ is the derivation on the ring ${\cal A}'$ such
 that
 $$\delta'(h\otimes 1)=1\otimes dd^ch,\,\,\delta'(1\otimes dd^ch)=0$$
 for any $h\in A$.
\item
   $\pi\colon{\cal A}'\to{\cal A}$ is the ring homomorphism such that
   $$h_i\otimes 1\mapsto h_i,\,\,\,1\otimes dd^ch_i\mapsto dd^ch_i$$
 \end{enumerate}
  Theorem \ref{thmMonge} follows from proposition \ref{uttPi}.
  Proposition \ref{uttPi} says that
  the derivation $\delta'$
  survives on the quotient ring ${\cal A}$ of the ring ${\cal A}'$.
 \begin{utt}
  Let  $\pi(F)=0$;
  then $\pi\delta'(F)=0$.
   \label{uttPi}
 \end{utt}
 Proposition \ref{uttPi} follows from (see below) proposition \ref{uttCurrent}.
  \begin{thm}
  {\rm(1)} For any $\tau\in{\cal A}$
  there exists a unique {\rm PPH}-cycle $X=\iota(\tau)$ such that $\tau=\bar X$.

  {\rm(2)}
  If $\tau\in{\cal T},\:\deg\tau=2k$
  then $\iota(\tau)$ is a $(2n-k)$-dimensional {\rm PPH}-cycle of degree $0$.
 \label{thmCurrent}
 \end{thm}
\proof
The uniqueness of {\rm PPH}-cycle is obvious.

First we prove the existence of $\iota(\tau)$ for
$\tau\in{\cal T}$ and the assertion (2).
The proof is by induction on $\deg(\tau)$.

If $\deg(\tau)=0$ then $\tau$ is a constant function $f(x)=c$
and $\iota(\tau)$ is the $2n$-cycle $X_\Delta=c$.

 Let $\deg(\nu)=2k-2$ and $\tau= dd^ch\wedge\nu$,
 where $h\in A$.
 By the inductive assumption,
 $\iota(\nu)$ is a $(2n-k+1)$-dimensional {\rm PPH}-cycle of degree $0$.

 By proposition \ref{uttEst2},
 it follows that $\tau=\overline{D_c(h\iota(\nu))}$.
 The degree of $(2n-k)$-dimensional PPH-cycle $D_c(h\iota(\nu))$ is $0$.
 So we can put $\iota(\tau)=D_c(h\iota(\nu))$.
 The theorem for $\tau\in{\cal T}$ is proved.

 The $\cal A$ as $\cal B$-module is generated by
 the elements of the ring ${\cal T}$.
 Similarly the $\cal B$-module of PPH-cycles is generated by
 PPH-cycles of degree $0$.
 So the map $\iota$ can be continued
 as the homomorphism of $\cal B$-modules.
  \begin{cor}
  If $\tau\in{\cal T}$ and $h\in A$ then
   $\iota(h\tau)=D_c(h\iota(\tau))$
 \label{corCurrent}
 \end{cor}
\begin{utt}
   $\iota\pi\delta'=D_c\iota\pi$
 \label{uttCurrent}
 \end{utt}
\proof
On elements of the ring ${\cal T}'$ both parts of the required equality are zeroes.
Let $\Upsilon\in{\cal T}'$ and $H\in{\cal B}'$ be an  element of the first degree.
If $\pi(1\otimes\Upsilon)=\nu$ and $\pi(H\otimes 1)=h$,
 then $\nu\in{\cal T}$ and $h\in A$.
  Any element of the ring ${\cal A}'$ is a linear combination of elements of the form
 $H^k\otimes\Upsilon$ ($k$ is not fixed).
 So we must prove that  $\iota\pi\delta'(H^k\otimes\Upsilon)=D_c\iota(h^k\nu)$.

  Now using the notation from the beginning of subsection, the Leibniz product rule,
  and the corollaries \ref{corCurrent} and \ref{corhk}, we get
  \begin{multline*}
 \iota\pi\delta'(H^k\otimes\Upsilon)= \iota\pi\delta'\left((H\otimes 1)^k (1\otimes\Upsilon)\right)=
 \iota\pi\left(\delta'\left((H\otimes 1)^k\right)(1\otimes\Upsilon)\right)=\\
 \iota\pi\left(k(H\otimes 1)^{k-1}\delta'(H\otimes 1)(1\otimes\Upsilon)\right)
 =\iota\pi\left(k(H\otimes 1)^{k-1} (dd^ch\otimes1)(1\otimes\Upsilon)\right)=
  \\\iota\left(kh^{k-1}\nu dd^ch \right)=kh^{k-1}\iota(dd^c(h\nu))
  =kh^{k-1}D_c\left(h\iota(\nu)\right)=D_c\iota(h^k\nu).
 \end{multline*}
%
%
\selectlanguage{english}

\par\bigskip
\par\bigskip
\par\bigskip
\begin{flushleft}
{\small
Institute for Information Transmission Problems,
B.Karetny per. 19, 101447 Moscow, Russia\\
E-mail address: kazbori@iitp.ru
}
\end{flushleft}
\end{document}